\documentclass[12pt]{amsart}
  \usepackage{graphics}
\usepackage{amsthm}

\begin{document}

\newtheorem{lemma}{Lemma}
\newtheorem{prop}[lemma]{Proposition}
\newtheorem{cor}[lemma]{Corollary}
\newtheorem{thm}{Theorem}
\newtheorem*{thm*}{Theorem}

\theoremstyle{definition}
\newtheorem{rem}[lemma]{Remark}
\newtheorem{rems}[lemma]{Remarks}
\newtheorem{defi}[lemma]{Definition}
\newtheorem{ex}{Example}
\newtheorem{convention}[lemma]{Convention}

\newcommand{\ra}{\longrightarrow}
\newcommand{\hyp}{h}
\newcommand{\cH}{\hyp}
\newcommand{\oh}{\overline{h}}
\newcommand{\ok}{\overline{k}}
\newcommand{\St}{\hbox{st}}
\newcommand{\lk}{\hbox{lk}}
\title[Hilbert space compression]{
Hilbert space compression and exactness of discrete groups.}
\author{Sarah Campbell and Graham Niblo}
\email{S.J.Campbell@maths.soton.ac.uk,
G.A.Niblo@maths.soton.ac.uk}
\date{\today}
\thanks{The first author was supported by an EPSRC postgraduate studentship.}
\begin{abstract}We show that the Hilbert space compression of any finite dimensional  CAT(0) cube complex is $1$ and deduce that any discrete group acting
properly, co-compactly on a CAT(0) cube complex is exact. The class of groups covered by this theorem includes free groups, finitely generated Coxeter groups, finitely generated  right angled Artin groups, finitely presented groups satisfying the  B(4)-T(4) small cancellation condition and all those  word-hyperbolic groups satisfying the B(6) condition. Another family of examples is provided by certain canonical
surgeries defined by link diagrams.\end{abstract}
\maketitle
\section*{Introduction}\label{intro}

We say that a group $\Gamma$ is exact if the operation of taking the reduced crossed product with $\Gamma$ preserves exactness of short exact sequences of $\Gamma$-$C^*$-algebras. In other words, $\Gamma$ is exact if and only if for every exact sequence of $\Gamma$-$C^*$-algebras
$$
0 \ra B \ra C \ra D\ra 0
$$
the sequence
$$
0 \ra C^*_r(\Gamma, B) \ra C^*_r(\Gamma, C) \ra C^*_r(\Gamma , D)\ra 0
$$
of crossed product algebras is exact. Kirchberg and Wassermann \cite{KW} proved that when 
$\Gamma$ is discrete, it is exact 
if and only if its reduced $C^*$-algebra $C^*_r(\Gamma)$ is exact. This means that the functor $B \mapsto 
C^*_r(\Gamma) \otimes_{\text{min}} B$ is exact, i.e. preserves exactness of sequences of 
$C^*$-algebras. 

Exact groups satisfy the coarse Baum Connes conjecture and the property was first made prominent by the work  of Kirchberg and Wassermann \cite{KW}, and studied by several authors \cite{AD,GK,HR,O,Yu1,Yu2}. Examples of exact groups include amenable groups, and groups of finite asymptotic dimension.

Exact groups admit a uniform embedding into Hilbert space and in \cite{GK2} Guentner and Kaminker 
introduced a numerical quasi-isometry invariant of a finitely generated group, whose values parametrize the difference between the group being uniformly embeddable in a Hilbert space and the reduced C*-algebra of the group being exact.

\begin{thm*}[Guentner, Kaminker]
Let $G$ be a discrete group. If the Hilbert space compression of $G$ is strictly greater than $1/2$ then $G$ is exact.
\end{thm*}

We will define Hilbert space compression later, but note here that it is a measure  of the amount of distortion that is necessary when trying to embed the group in a Hilbert space via a large scale Lipschitz map. Guentner and Kaminker illustrated their theorem by showing that the Hilbert space compression of a finite rank free group  is $1$ thus giving a new proof of exactness for free groups. It should be noted that they did not construct an embedding of the free group in a Hilbert space with (asymptotic) compression $1$, but rather, thinking of the group as a tree via its Cayley graph, produced a family of large scale Lipschitz embeddings with asymptotic  compression arbitrarily close to $1$. Those familiar with CAT(0) complexes would recognize that the first of their embeddings (with asymptotic compression $1/2$) can be used without change to embed the vertex set of a CAT(0) cube complex into a Hilbert space with asymptotic compression $1/2$ though this is not in itself enough to establish exactness for a group acting on the cube complex. Guentner and Kaminker showed that in the case of a tree the embedding can be modified to obtain new embeddings with asymptotic compression arbitrarily close to $1$. 

The main purpose of this note is to show how to adapt the construction from \cite{GK2} to the class of finite dimensional CAT(0) cube complexes. In the case of a tree one uses the fact that there is a unique edge geodesic joining any two points in the tree; the same is of course not true for CAT(0) cube complexes of dimension at least 2 so the embedding and the argument need to be modified appropriately. In place of unique edge geodesics we will use the normal cube paths originally introduced in \cite{NR1} to establish biautomaticity for groups acting freely and properly discontinuously on CAT(0) cube complexes. 

\begin{thm*}[Theorem A] Let $X$ be a finite dimensional CAT(0) cube complex. The HIlbert space compression of $X$ is $1$.
\end{thm*}

 In \cite{GK} it is shown that Hilbert space compression is a quasi-isometry invariant so if a discrete group $G$ acts freely and co-compactly on a CAT(0) cube complex it follows that the group (regarded as a metric space via the word length metric) has Hilbert space compression $1$. Since $1>1/2$ we obtain:

\begin{thm*}[Theorem B] If $G$ is a group acting properly and co-compactly on a CAT(0) cube complex then $G$ is exact.
\end{thm*}

The paper is organised as follows: In section \ref{cubes} we recall the definition of a CAT(0) cube complex and, stating the definitions, show how to construct a large scale Lipschitz embedding of such a complex in an associated Hilbert space, with asymptotic compression $1/2$. In section \ref{normal} we outline some preliminary results concerning the existence and properties of normal cube paths in a CAT(0) cube complex. The results in this section are taken from \cite{NR1}. In section \ref{Lipschitz} we define a family of embeddings $\{f_\epsilon\mid 0<\epsilon<1/2\}$ of the vertices of a cube complex $X$ into the Hilbert space of square summable real valued functions on the set of hyperplanes of $X$. We also show that these embeddings are large-scale Lipschitz. In section \ref{compression} we show that the compression of each map $f_\epsilon$ is $1/2+\epsilon$ and deduce that the Hilbert space compression of the metric space $(X^{(0)}, d_1)$ is $1$, where $X^{(0)}$ denotes the vertex set of $X$ and $d_1$ is the edge metric. In section \ref{exact} we deduce the exactness of groups acting properly and co-compactly on a CAT(0) cube complex.

The class of groups covered by this theorem includes free groups, finitely generated Coxeter groups \cite{NR3}, and finitely generated  right angled Artin groups (for which the Salvetti complex is a CAT(0) cube complex). A rich class of interesting examples is furnished by Wise, \cite{W}, in which it is shown that many small cancellation groups  act properly and co-compactly on CAT(0) cube complexes. The examples include  every finitely presented group satisfying the  B(4)-T(4) small cancellation condition and all those  word-hyperbolic groups satisfying the B(6) condition. Finally many 3-manifolds admit  decompositions as CAT(0) cube complexes, so their fundamental groups are also covered by the theorem, a family of examples is provided by certain canonical
surgeries defined by link diagrams (see \cite{AR} and \cite{AR2}). Classical examples are furnished by groups acting simply transitively on buildings with the structure of a product of trees.

The authors wish to thank Jacek Brodzki and Claire Vatcher for many interesting and illuminating conversations during the course of this research.

\section{CAT(0) cube complexes}\label{cubes}
A {\em cube complex} $X$ is a metric polyhedral complex in which
each cell is isometric to the Euclidean cube $[-1/2,1/2]^n$, and
the gluing maps are isometries. If there is a bound on the
dimension of the cubes then the complex carries a complete
geodesic metric, \cite{BH}.

\begin{rem} A cube complex is {\em non-positively curved} if for
any cube $C$ the following conditions on the link of $C$,
$\lk{C}$, are satisfied:
\begin{enumerate}
\item({\em no bigons}) For each pair of vertices in
$\lk{C}$ there is at most one edge containing them.
\item ({\em no triangles}) Every
edge cycle of length three in $\lk{C}$ is contained in a 2-simplex
of $\lk{C}$.
\end{enumerate}
\end{rem}

The following theorem of Gromov relates the combinatorics and the
geometry of the complex.

\begin{lemma} {\rm (Gromov, \cite{G})}
A cube complex $X$ is locally $CAT(0)$ if and only if it is
non-positively curved, and it is CAT(0) if and only if it is
non-positively curved and simply connected.
\end{lemma}

\begin{ex}
Any graph may be regarded as a 1-dimensional cube complex, and the
curvature conditions on the links are trivially satisfied. The
graph is CAT(0) if and only if it is a tree. Euclidean space also
has the structure of a CAT(0) cube complex with its vertices at
the integer lattice points.
\end{ex}

A {\em midplane} of a cube $[-1/2,1/2]^n$ is its intersection with
a codimension 1 coordinate hyperplane. So every $n$-cube contains
$n$ midplanes each of which is an $(n-1)$-cube, and any $m$ of
which intersect in a $(n-m)$-cube. Given an edge in a
non-positively curved cube complex, there is a unique codimension
1 hyperplane in the complex which cuts the edge transversely in
its midpoint. This is obtained by developing the midplanes in the cubes containing the edge. In the case of a
tree the hyperplane is the midpoint of the edge, and in the case
of Euclidean space it is a geometric (codimension-1) hyperplane.

In  general a hyperplane is analogous to an immersed codimension
1 submanifold in a  Riemannian manifold and in a CAT(0) cube complex one can show that the immersion is
 a local isometry. An application of the Cartan-Hademard
theorem, then shows that the hyperplane is
isometrically embedded. Furthermore any hyperplane in a CAT(0)
cube complex separates it into two components referred to as the
half spaces associated with the hyperplane. This is a consequence
of the fact that the complex is simply connected. The hyperplane
gives rise to 1-cocycle which is necessarily trivial, and hence
the hyperplane separates the space.

The set of vertices of a CAT(0) cube complex $X$ can be can be viewed as a
discrete metric space, where the metric $d_1(u,v)$ is given by the
length of a shortest edge path between the vertices $u$ and $v$. We will refer to this as the $\ell^1$ metric on the vertices. Alternatively we can measure the distance by restricting the path metric on $X$ to obtain the $\ell^2$ metric on the vertices. If $X$ is finite dimensional these metrics
are
quasi-isometric, and we have $d(u,v)\leq d_1(u,v)\leq \sqrt{n}d(u,v)$
where $d$ denotes the CAT(0) (geodesic) metric on $X$ and $n$ is the dimension of the
complex. 

Sageev \cite{S} observed that the shortest path in the 1-skeleton
crosses any hyperplane at most once, and since every edge crosses
exactly one hyperplane, the $\ell^1$ distance between two vertices is the
number of hyperplanes separating them. 

Finally we will need the concept of a median. In any CAT(0) cube complex there is a well defined notion of an interval; given any two vertices $u,v$ the interval between them, denoted $[u,v]$ consists of all the vertices which lie on an edge geodesic from $u$ to $v$. Given any three vertices $u,v,w$ there are three intervals $[u,v], [v,w], [w,u$ and the intersection of these three intervals is always a single point $m$ known as the median of the triple$u,v,w$ (see \cite{R} for details). It has the following important property: If we consider the hyperplanes which separate the pair $u,v$ and those which separate the pair $u,w$ the intersection of these two families consists of precisely the hyperplanes which separate $u$ and the median $m$. Furthermore the hyperplanes which separate $v$ from $w$ are precisely those hyperplanes which separate $m$ from $v$ together with those which separate $m$ from $w$ so we have $d_1(v,w)=d_1(v,m)+d_1(m,w)=d_1(v,u)+d_1(w,u)-2d_1(m,u)$. We will use this fact in section \ref{compression}.


In \cite{NRo} it was shown how to use the hyperplane structure of a CAT(0) cube complex $X$ to obtain an $\ell^1$ embedding of the cube complex in the Hilbert space $\ell^2(H,\mathbb R)$ of square summable (real valued)  functions on the set  $H$ of hyperplanes in $X$. An alternative description of the embedding, based on the one used in \cite{GK2} in the context of a tree, is as follows:

Choose a basepoint $v$ in $X^{(0)}$ and for each vertex $w\in X^{(0)}$ set $H_w=\{h\in H\mid h\hbox{ separates $v$ and $w$ }\}$. Define $f_w: H\longrightarrow \mathbb R$ by $f_w=\sum\limits_{h\in H_w}\delta_h$ where $\delta_h$ denotes the characteristic function of the singleton $\{h\}\subset H$.

It is easy to see that the function $f_w$ is $\ell^1$ and therefore $\ell^2$ and since the Hilbert space is contractible (in fact uniquely geodesic) the map extends to an embedding of $X$ in $\ell^2(H,\mathbb R)$. If $X$ is a cube then this embedding is isometric, however in the case of a tree (consisting of more than a single edge) then it is not.

\begin{ex}
Let $T$ be the tree consisting of two edges $e_s, e_t$ both adjacent to a vertex $v$, and with the other two vertices labelled $s,t$. The tree has two hyperplanes, corresponding to the midpoints of the two edges, so that $\ell^2(T,\mathbb R)\sim \mathbb Re_s\bigoplus\mathbb Re_t$. The vertex $v$ is not separated from itself by either of the hyperplanes so we have $f_v=0$. The vertex $s$ is only separated from $v$ by the hyperplane $s$ so we have $f_s=\delta_{e_s}$ and similarly $f_t=\delta_{e_t}$. Now in the tree we have $d_1(s,t)=d_2(s,t)=2$ however in the Hilbert space we have $d_1(f_s, f_t)=2\not=\sqrt{2}=d(f_s, f_t)$, where we have used $d_1$ to denote the $\ell^1$ metric and $d$ to denote the Hilbert metric. 
\end{ex}


Although the embedding defined above is not necessarily an isometry it is relatively easy to show that it is a large scale Lipshcitz map, and we can measure the distortion of such a map in terms of its compression:

\begin{defi} A function $f:X\rightarrow Y$ is large-scale
Lipschitz if there exist $C>0$ and $D\geq 0$ such that $d_Y(f(x),f(y))\leq Cd_X(x,y)+D$
Following Gromov, the compression $\rho(f)$
of $f$ is given by $\rho_f(r)=\inf_{d_X(x,y)\geq
r}d_Y(f(x),f(y))$.  Assuming that $X$ is unbounded the asymptotic compression $R_f$ is  given by $$R_f=\liminf_{r\rightarrow\infty}
\frac{\log\rho^*_f(r)}{\log r}$$ where
$\rho^*_f(r)=\max\{\rho_f(r),1\}.$
\end{defi}

In the case of the embedding of the vertices described above the map is large scale Lipschitz with  $C=1, D=0$. The argument used by Guentner and Kaminker \cite{GK2} to compute the asymptotic compression of the embedding of a tree goes through without change to our more general context to show that the asymptotic compression is $1/2$. (It should be noted here that we are regarding the cube complex as a metric space via the $\ell^1$ metric not the (geodesic) $\ell^2$ metric.)

In order to obtain large scale Lipschitz embeddings with asymptotic compression close to $1$ we need to  adapt the embedding described above. The idea, taken from \cite{GK2} is to weight the functions $\delta_h$ according to how far the hyperplane is from the basepoint. Whereas in the case of a tree the hyperplanes which separate two vertices are linearly ordered in  a higher dimensional cube complex they are not and there are several partial orders one could use in modifying the argument. It turns out that the appropriate ordering is furnished by the normal cube paths introduced in \cite{NR1} and we describe these next.


\section*{Normal cube paths}\label{normal}



\begin{defi} A cube path is a sequence of cubes $\mathcal C=\{C_0, \ldots C_n\}$, each of
dimension at least 1, such that each cube meets its successor in a
single vertex, $v_i=C_{i-1}\cap C_i$ and such that for $1\leq
i\leq n-1$, $C_i$ is the (unique) cube of minimal dimension
containing $v_i$ and $v_{i+1}$.  Note that $v_i$ and $v_{i+1}$ are
diagonally opposite vertices of $C_i$. We define $v_0$ to be the
vertex of $C_0$ which is diagonally opposite $v_1$, and $v_{n}$ to
be the vertex of $C_n$ diagonally opposite $v_{n-1}$. We call the
$v_i$, vertices of the cube-path, with $v_0$ the initial vertex
and $v_{n}$ the terminal vertex. Given a cube path from $u$ to $v$
we can construct edge paths from $u$ to $v$ which travel via the
edges of the cubes $C_i$ so every hyperplane separating $u$ from
$v$ must intersect at least one of the cubes $C_i$. We say the
cube path is normal if $C_{i+1}\cap\St(C_{i})=v_i$ for each $i$,
where $\St(C_i)$ is the union of all cubes which contain $C_i$ as
a face (including $C_i$ itself).
\end{defi}

In \cite{NR1} it was shown that given any two vertices
$u,v$ there is a unique normal cube path $\mathcal C=\{C_0,\ldots, C_n\}$ from $u$ to $v$. We will need the following key facts about normal cube paths all of which may be found in \cite{NR1}.

\begin{lemma} \label{technical} Let $s,t,v_0$ be vertices of a CAT(0) cube complex with $s$ and $t$ diagonally adjacent across some cube $E_0$. Let
$s=s_0, s_1, \ldots, s_m=v$, $t=t_0, t_1, \ldots, t_n=u$ be the
vertices of the (unique) normal cube paths from $s$ to $v_0$ and
from $t$ to $v_0$ respectively. Let $\{C_i\mid i=1,\ldots m\}$ be
the cubes on the normal cube path from $s$ to $v_0$ and $\{D_j\mid
j=1,\ldots n\}$ be the cubes on the normal cube path from $t$ to
$v_0$. Then:
\begin{enumerate}
\item  Each hyperplane separating $s$ from $v_0$ intersects exactly one of the cubes  $C_i$ and each hyperplane separating $t$ from $v_0$ intersects exactly one of the cubes $D_j$.
\item For each $i\leq \min\{m,n\}$ there is a cube $E_i$ such that $s_i$ is diagonally adjacent to $t_i$ across $E_i$ .
\end{enumerate}
\end{lemma}

We will need the following technical lemma:

\begin{lemma}\label{controllingweights}  Let $s, t, v$ be vertices of the CAT(0) cube complex $X$
with $s$ and $t$ diagonally opposite across some cube $E_0$. Let
$s=s_0, s_1, \ldots, s_m=v$, $t=t_0, t_1, \ldots, t_n=v$ be the
vertices of the (unique) normal cube paths from $s$ to $v$ and
from $t$ to $v$ respectively. Let $\{C_i\mid i=1,\ldots m\}$ be
the cubes on the normal cube path from $s$ to $v$ and $\{D_j\mid
j=1,\ldots n\}$ be the cubes on the normal cube path from $t$ to
$v$. If $h$ is a hyperplane in $X$ which separates both $s$ and
$t$ from $v$ and which intersects the cube $C_i$ then $h$ also
intersects one of the cubes $D_{i-1}, D_i, D_{i+1}$.
\end{lemma}

\begin{proof}  By lemma \ref{technical} the hyperplane $h$ can only (and must) intersect the normal cube path from $t$ to $v$  in one
of the cubes $D_j$, and the hypothesis that $s=s_0$ and $t=t_0$
are diagonally opposite across the cube $E_0$ ensures that for each $i\leq \min\{m,n\}$ $s_i$ is
diagonally opposite to $t_i$ across some cube $E_i$.

Now $h$ separates $s_{i-1}, s_i$ and also separates $t_{j-1}, t_j$. let  $k=\min\{i,j\}$. Assume first
that $h$ separates $s_{k-1}$ and $s_k$ so $i=k\leq j$; if $h$ also
separates $t_{k-1}$ and $t_{k}$ then $h$ crosses $D_{k}=D_i$ as
required. If on the other hand $h$ does not separate $t_{k-1}$
from $t_{k}$ then, since it does not separate $t$ from $t_{k-1}$
by the minimality of $k$, but does separate $t$ from $v$, it must
also separate $t_{k}$ from $v$. Now we construct an edge path from
$t_{k}$ to $v$ as follows. First cross over the cube $E_{k}$ from
$s_{k}$ to $t_{k}$ then follow a path through the cubes $C_{k+1},
C_{k+2},\ldots, C_m$ to $v$. This gives an edge path from $t_{k}$
to $v$ so it must cross $h$. However none of the cubes
$C_{k+1},\ldots, C_m$ intersect $h$ so $h$ must cross $E_{k}$ and
hence $h$ is adjacent to $t_{k}$. But as $h$ separates $t_{k}$
from $v$ and is adjacent to $t_{k}$ it must cross the first cube
($D_{k+1}$) on the normal cube path from $t_{k+1}$ to $v$ as
required. The case when $h$ separates $t_{k-1}$ and $t_k$ so $i=k$
but does not separate $s_{k-1}$ and $s_{k}$ is argued in exactly
the same way reversing the roles of $s$ and $t$, $C$ and $D$ and
so on.
\end{proof}

From now on we fix a vertex $v$ as a basepoint and for each vertex $s$ we define an integer-valued weight function $w_s$ on the set of hyperplanes as follows. Let $\mathcal C=\{C_0,\ldots, C_n\}$ be the unique normal cube path from $s$ to $v$. If the hyperplane $\hyp$ separates $s$ and $v$ then set  $w_s(\hyp)=i+1$ where  $h$ intersects the cube $C_i$, otherwise set $w(\hyp)=0$. Hence $w_s$ has finite support. From Lemma \ref{controllingweights} we get:

\begin{cor}\label{controllingweightscor} If $s$ and $t$ are adjacent in $X$ and $h$ separates both
$s$ and $t$ from $v$ then $|w_t(h)-w_s(h)|\leq1$.
\end{cor}

\begin{proof} If the normal cube path from $s$ to $v$ is denoted by the
cubes $C_i$ as above and the normal cube path from $t$ to $v$ is
denoted by $D_j$ then $h$ intersects precisely the cubes
$C_{w_s(h)}$ and $D_{w_t(h)}$ so by the lemma  
$D_{w_t(h)}=D_{w_s(h)\pm 1}$, and $w_t(h)=w_s(h)$ or $w_t(h)=w_s(h)\pm 1$ as
required.
\end{proof}

Note that in the statement of the corollary ``adjacent'' may be taken to mean adjacent across the diagonal of any cube, however in our application we will only need it to mean that $s$ and $t$ are vertices of a common edge.

\section{The large scale Lipschitz embeddings}\label{Lipschitz}

As in the last section we fix a CAT(0) cube complex $X$ (not necessarily
finite dimensional) and a base vertex $v$. We will show how to construct  a family (indexed by the interval $(0, 1/2)$) of large scale Lipschitz embeddings of the vertex set $X^{(0)}$ into the Hilbert space of $\ell^2$ functions from the set $H$ of hyperplanes in $X$ to $\mathbb R$.

%

For each $\epsilon\in(0,1/2)$ define  $f_\epsilon(s)=\sum\limits_{h\in {\hyp}}{w_s(h)}^\epsilon\delta_{h}$. As noted before since the $s$-weight of a hyperplane is $0$ unless the hyperplane is one of the finitely many separating $s$ from the basepoint $v$, this sum is always finite and therefore is an element of $\ell^2(H,\mathbb R)$.


 In order to show that $f_\epsilon$ is a large scale Lipschitz map it suffices to show that there is a constant $C$ such that whenever $d_1(s,t)=1$, $\parallel
f_\epsilon(s)-f_\epsilon(t)\parallel^2\leq C$.

\begin{lemma}
For each $\epsilon\in (0,1/2)$ there is a constant $C$ such that for any vertices $s,t\in X^{(0)}$ with $d_1(s,t)=1$ we have $\parallel
f_\epsilon(s)-f_\epsilon(t)\parallel^2\leq C$.
\end{lemma}

\begin{proof}
Let $\hyp_0$ be the hyperplane cutting the edge joining $s,t$. Assume, without loss of generality that $h$ separates $t$ from $v$ but not $s$ from $v$ so that $d_1(s,v)+1=d_1(t,v)$ and the set of hyperplanes separating $t$ from $v$ is the union of the set $\{\hyp_1,\ldots,\hyp_m\}$ of the hyperplanes separating $s$ from $v$ together with $\hyp$.

We need to compute

\begin{equation}
\parallel f_\epsilon(s)-f_\epsilon(t)\parallel^2=\sum\limits_{i=0}^m[w_s(\hyp_i)^\epsilon-w_t(\hyp_i)^\epsilon]^2=1^{2\epsilon}+\sum\limits_{i=1}^m[w_s(\hyp_i)^\epsilon-w_t(\hyp_i)^\epsilon]^2
\end{equation}

Now according to corollary  \ref{controllingweightscor}  we have $|w_t(h_i)-w_s(h_i)|\leq1$. Suppose that for a particular hyperplane $\hyp_i$ we have $w_s(\hyp_i)=k$  so that $w_t(\hyp_i)$ takes one of the values $k-1, k, k+1$ and $[w_s(\hyp_i)^\epsilon-w_t(\hyp_i)^\epsilon]^2$ takes one of the values $[k^\epsilon-(k+1)^\epsilon]^2, [k^\epsilon-k^\epsilon]^2, [k^\epsilon-(k-1)^\epsilon]^2$

An elementary calculation of the first derivative shows that the function $X\mapsto [X^\epsilon-(X+1)^\epsilon]^2$.  is strictly increasing so  we have $[k^\epsilon-(k+1)^\epsilon]^2>[(k-1)^\epsilon-k^\epsilon]^2=[k^\epsilon-(k-1)^\epsilon]^2>0=[k^\epsilon-k^\epsilon]^2$ hence  we have  $[w_s(\hyp_i)^\epsilon-w_t(\hyp_i)^\epsilon]^2\leq [w_s(\hyp_i)^\epsilon-(w_s(\hyp_i)+1)^\epsilon]^2$ and so 

\begin{equation}
\sum\limits_{i=1}^m[w_s(\hyp_i)^\epsilon-w_t(\hyp_i)^\epsilon]^2\leq \sum\limits_{i=1}^m  [w_s(\hyp_i)^\epsilon-(w_s(\hyp_i)+1)^\epsilon]^2.
\end{equation}

We can split the final sum as a double sum taken over all hyperplanes with a given $s$-weight. Let $J$ denote the set of all $s$-weights.

\begin{equation}
\sum\limits_{i=1}^m  [w_s(\hyp_i)^\epsilon-(w_s(\hyp_i)+1)^\epsilon]^2=
\sum\limits_{j\in J}\sum\limits_{w_s(\hyp_i)=j}  [w_s(\hyp_i)^\epsilon-(w_s(\hyp_i)+1)^\epsilon]^2=
\sum\limits_{j\in J}\sum\limits_{w_s(\hyp_i)=j}  [j^\epsilon-(j+1)^\epsilon]^2.
\end{equation}

Since the cube complex has dimension $n$  we can cross at most $n$ hyperplanes in any given cube so the number of hyperplanes with $w_s(\hyp)=j$ is at most $n$ for any $j$ and

\[
\sum\limits_{j}\sum\limits_{w_s(\hyp_i)=j}  [j^\epsilon-(j+1)^\epsilon]^2
\leq \sum\limits_{j}n[j^\epsilon-(j+1)^\epsilon]^2.
\]

Putting $w_j=j^\epsilon$  and adding additional positive terms we see  that 

\[
\parallel f_\epsilon(s)-f_\epsilon(t)\parallel^2\leq n\sum\limits_{j=0}^\infty [w_j-w_{j+1}]^2.
\]

As remarked  in \cite{GK2} the series  $\sum\limits_{j=0}^\infty [w_j-w_{j+1}]^2$ converges  so that putting $C=n\sum\limits_{j=0}^\infty [w_j-w_{j+1}]^2$ we get $\parallel f_\epsilon(s)-f_\epsilon(t)\parallel^2\leq C$ and  $f_\epsilon$ is
large scale Lipschitz as required.
\end{proof}

\section{Hilbert space compression}\label{compression}

While establishing that the map is large scale Lipschitz required us to show that $\parallel
f_\epsilon(s)-f_\epsilon(t)\parallel^2$ is small for nearby vertices, to establish that the embedding has  large asymptotic compression requires us to show that $\parallel
f_\epsilon(s)-f_\epsilon(t)\parallel^2$ is relatively large for points $s,t$ which are sufficiently far apart.

Specifically we will prove:

\begin{lemma}
For any positive $r$ and any $\epsilon \in (0,1/2)$ there is a constant $C_\epsilon$ such that $\parallel
f_\epsilon(s)-f_\epsilon(t)\parallel^2\geq C_\epsilon
r^{1+2\epsilon}$. Hence $\rho_{f_\epsilon}(r)\geq
\sqrt{C_\epsilon}r^{1/2+\epsilon}$

\end{lemma}

\begin{proof}
Let $D=d_1(s,t)\geq r$ and assume $d(1,s)\leq d(1,t)$ so that, letting $m$ denote the median of the triple $1,s,t$, we have $d(m,t)\geq d(m,s)$. It follows that $d(m,t)\geq \sharp(\frac{D}{2})\geq \sharp(\frac{r}{2})$ where $\sharp(n)$ denotes the smallest integer greater than
$n$. Hence there are at least $\sharp(\frac{r}{2})$ hyperplanes which separate $t$ from $1$ but which do not separate $s$ from $1$. We will denote these hyperplanes $h_1, h_2, \ldots, h_{\sharp(\frac{r}{2})}$.
 Now consider the normal cube path $C_0, C_1,\ldots C_n$ from $t$ to $1$. As noted in lemma \ref{technical} each of the hyperplanes $h_i$ must intersect exactly one of the cubes $C_j$, and by definition $w_t(h_i)=(j+1)$. By relabelling if necessary we may assume that the $t$-weight increases (not necessarily strictly) with the index $i$ of the hyperplane, and given that the cube complex has dimension $n$ at most $n$ of the hyperplanes can have the same $t$-weight, i.e., at most $n$ of the hyperplanes have weight $1^\epsilon$ and the others have weight at least $2^\epsilon$; at most $n$ of the remaining hyperplanes can have weight $2^\epsilon$ and the others have to have weight at least $3^\epsilon$ and so on. Recall that for each of these hyperplanes $w_s(h_i)=0$ by hypothesis so, writing $\sharp(\frac{r}{2})=kn+m$ for some integer $0\leq m<n$ we have
 
 $$\parallel
f_\epsilon(s)-f_\epsilon(t)\parallel^2\geq w_t{\pi(h_1)}^2+\ldots+w_t{\pi(h_{\sharp(\frac{r}{2})})}^2
\geq n(1^{2\epsilon}+2^{2\epsilon}+\ldots+k^{2\epsilon})+m(k+1)^{2\epsilon}.$$

We will now show that the RHS of this equation is greater than the
expression
$$\frac{1}{n}(w_1^2+w_2^2+\ldots+w_{\sharp(\frac{r}{2})}^
2)$$
$$=\frac{1}{n}(w_1^2+\ldots+w_n^2+w_{n+1}^2+\ldots+w
_{2n}^2+w_{2n+1}^2+\ldots
+w_{kn}^2+w_{kn+1}^2+\ldots+w_{kn+m=\sharp(\frac{r}{2})}^
2)$$

\textbf{Claim:} For any $i\geq1$,

\begin{displaymath}
ni^{2\epsilon}>\frac{1}{n}[(i-1)n+1)^{2\epsilon}+\ldots+(in)^{2\epsilon}]
\end{displaymath}

Since $\epsilon<\frac{1}{2}$ and $n\geq 1$ we have $ni^{2\epsilon}>n^{2\epsilon}i^{2\epsilon}=(in)^{2\epsilon}$

On the other hand, since $\epsilon>0$ and $in>ik$ for all $k<n$ we have
$$
\frac{1}{n}[((i
-1)n+1)^{2\epsilon}+\ldots+(in)^{2\epsilon}]
<\frac{1}{n}(n(in)^{2\epsilon})=(in)^{2\epsilon}
$$

So $$ni^{2\epsilon}>(in)^{2\epsilon}>\frac{1}{n}[{((i
-1)n+1)}^{2\epsilon}+\ldots+{in}^{2\epsilon}]$$

\textbf{Claim:}
\begin{displaymath}
mw^2_{k+1}>\frac{1}{n}\left(w^2_{kn+1}+\ldots+w^2_{kn+m
}\right)
\end{displaymath}

We have that  $$mw^2_{k+1}=m(k+1)^{2\epsilon}\geq
m^{2\epsilon}(k+1)^{2\epsilon}=(mk+m)^{2\epsilon}$$

Looking at the RHS of the statement of the claim we have:
\begin{eqnarray*}
\frac{1}{n}\left(w^2_{kn+1}+\ldots+w^2_{kn+m}\right)&=&\frac{1
}{n}\left((kn+1)^{2\epsilon}+\ldots+(kn+m)^{2\epsilon} \right)\\
&<&\frac{m}{n}(kn+m)^{2\epsilon}\hspace{5mm}{\rm
(since\,(kn+m)\,is\,the\,biggest\,term)}\\
&<&\left(\frac{m}{n}\right)^{2\epsilon}(kn+m)^{2\epsilon}
\hspace{5mm}({\rm since } \frac{m}{n}<1)\\
&=&(mk+\frac{m^2}{n})^{2\epsilon}\\ &<&(mk+m)^{2\epsilon}
\hspace{5mm}({\rm since } \frac{m}{n}<1)\\
\end{eqnarray*}
And so $$mw^2_{k+1}=m(k+1)^{2\epsilon}\geq
(mk+m)^{2\epsilon}>\frac{1}{n}\left(w^2_{kn+1}+\ldots+w^2_{kn+
m}\right)$$

Putting both claims together, we have that:
$$nw_1^2>\frac{1}{n}(w_1^2+\ldots+w_n^2)$$
$$nw_2^2>\frac{1}{n}(w_{n+1}^2+\ldots+w_{2n}^2)$$
\hspace{50mm}\vdots
$$nw_k^2>\frac{1}{n}(w_{(k-1)n+1}^2+\ldots+w_{kn}^2)$$
$$mw^2_{k+1}\geq\frac{1}{n}\left(w^2_{kn+1}+\ldots+w^2_{k
n+m}\right)$$

And so
$$nw_1^2+nw_2^2+\ldots+nw_k^2+mw_{k+1}^2>\frac{1}{n}
(w_1^2+w_2^2+\ldots+w^2_{\sharp(\frac{r}{2})})$$

Hence, \begin{eqnarray*}
\parallel f_\epsilon(s)-f_\epsilon(t)\parallel^2
&\geq&w_{\pi(h_1,t)}^2+\ldots+w_{\pi(h_{\sharp(\frac{r}{2}),t}
}^2\\ &\geq&
nw_1^2+nw_2^2+\ldots+nw_k^2+mw_{k+1}^2\\
&\geq&\frac{1}{n}(w_1^2+w_2^2+\ldots+w^2_{\sharp(\frac{r}
{2})})\\
&=&\frac{1}{n}(c_{\epsilon,1}^2+c_{\epsilon,2}^2+\ldots+c_{\epsilon,\sharp(\frac{r}{2})}^2)\\
&\geq& \frac{r^{2\epsilon+1}}{n(2^{2\epsilon+1})(2\epsilon+1)}\\
\end{eqnarray*}

\end{proof}

Now we obtain:

\begin{lemma}
For each $\epsilon$ the asymptotic compression of the map $f_\epsilon$ is at least $1/2+\epsilon$.
\end{lemma}

\begin{proof}We have
$$R_{f_\epsilon}=\liminf_{r\rightarrow\infty}\frac{\log
\rho_{f_\epsilon}(r)}{\log
r}\geq\liminf_{r\rightarrow\infty}\frac{\log
\sqrt{C_\epsilon}r^{1/2+\epsilon}}{\log r}=1/2+\epsilon$$
\end{proof}

\section{Exactness for groups acting properly and co-compactly on a CAT(0) cube complex}\label{exact}

The Hilbert space compression of a metric space is defined to be the supremum of the asymptotic compression of all possible large scale Lipschitz maps from the metric space to a Hilbert space so putting together the results of sections \ref{Lipschitz} and \ref{compression} we get

\begin{thm}[Theorem A] The Hilbert space compression of a finite dimensional CAT(0) cube complex $(X, d)$ is $1$
\end{thm}

\begin{proof} For each $\epsilon\in (0,1/2)$ we have constructed a large scale Lipschitz embedding $f_\epsilon$ of the metric space $(X^{(0)}, d_1)$ into the Hilbert space $\ell^2(H,\mathbb R)$ with compression at least $1/2+\epsilon$. Hence the Hilbert space compression of $(X^{(0)}, d_1)$ is $1$. Since $X$ is finite dimensional, of dimension $n$ say, we have $d(s,t)\leq d_1(s,t)\leq \sqrt{n}d(s,t)$ so $(X^{(0)}, d_1)$ is quasi-isometric to $(X,d)$, and since Hilbert space compression is a quasi-isometry invariant we obtain the result.
\end{proof}

Now suppose $G$ s a group acting properly and co-compactly on a CAT(0) cube complex $X$. Choose a finite generating set for $G$ and regard $G$ as a metric space via the edge metric on the Cayley graph. Then $G$ is quasi-isometric to $(X,d)$. According to Guentner and Kaminker HIlbert space compression is a quasi-isometry invariant so we obtain

\begin{cor} Let $G$ be a finitely generated group regarded as a metric space via the word metric with respect to some finite generating set. If $G$ acts properly and co-compactly on a CAT(0) cube complex then $G$  has Hilbert space compression $1$.
\end{cor}

Finally since Guentner and Kaminker showed that a discrete group with Hilbert space compression strictly greater than $1/2$ is exact we obtain:

\begin{thm*}[Theorem B] If $G$ is a group acting properly and co-compactly on a CAT(0) cube complex then $G$ is exact.
\end{thm*}



\begin{thebibliography}{99}



\bibitem{AD} C.~Anantharaman-Delaroche, J.~Renault, 
Monographies de L'Enseignement MathŽmatique, 36. L'Enseignement MathŽmatique, Geneva, 2000. 


\bibitem {AR}I.~R.~Aitchison, J.~H.~Rubinstein, An introduction to polyhedral metrics of non-positive curvature on 3-manifolds in "Geometry of Low-dimensional Manifolds: 2 Symplectic Manifolds and Jones-Witten Theory", London Mathematical Society Lecture Notes Series 151, Cambridge University Press (1990), pp 127--161.

\bibitem{AR2} I.~R.~Aitchison, J.~H.~ Rubinstein, Canonical surgery on alternating link diagrams in Knots 90, (1992) pp 543--558.

\bibitem{BH}
M.R. Bridson and A.~Haefliger, Metric spaces of non-positive curvature, {Grundlehren der mathematischen Wissenschaften Volume 319, 
Springer, 1999}

\bibitem{G} M.~Gromov, Hyperbolic groups in S.M. Gersten (ed.)
Essays on Group Theory
(Mathematical Sciences Research Institute Publications) Springer-Verlag (1987).

\bibitem{GK} E.~Guentner, J.~Kaminker, Exactness and the Novikov conjecture. Topology 41 (2002), no. 2, 411--418;  Addendum: Topology 41 (2002), no. 2, 419--420

\bibitem{GK2}E.~Guentner, J.~Kaminker, Exactness and Uniform Embeddability of Discrete Groups,
ArXiv preprint math.OA/0309166.

\bibitem{HR} N.~Higson, J.~Roe,  Amenable group actions and the Novikov conjecture. J. Reine Angew. Math. 519 (2000), 143--153.


\bibitem{KW}Kirchberg, Wasserman, Exact Groups and Continuous
Bundles of $C^\ast$-algebras, Mathematische Annalen, 1999.

\bibitem{NRo} G.~Niblo, M.~Roller, Groups acting on cubes and Kazhdan's property (T).
Proc. Amer. Math. Soc.  {\bf 126}  (1998),  no. 3, 693--699.

\bibitem{NR1}G.A~.Niblo, L.D.~Reeves. The geometry of cube
complexes and the complexity of their fundamental groups,Topology,
Vol. 37,No 3 (1998) pp621-633.


\bibitem{NR2}G.A.Niblo, L.D.Reeves. Groups acting on CAT(0) cube complexes, Geometry and Topology, Vol. 1 (1997) Paper no. 1, pages 1-7. 

\bibitem{NR3}G.A.Niblo, L.D.Reeves. Coxeter groups act on CAT(0) cube complexes,  Journal of Group Theory, 6, (2003), pp 309-413.


\bibitem{O} N.~Ozawa, Amenable actions and exactness for discrete groups, \emph{Comptes Rendus Acad. Sci. Paris} {\bf 330} (2000), 691--695. 

\bibitem{R} M.~A.~Roller,
Poc Sets, Median Algebras and Group Actions. An extended study of Dunwoody's construction and Sageev's theorem
{\em Habilitationeschrift, Regensberg\/} (1998)

\bibitem{S} M.~Sageev, Ends of group pairs and non-positively curved cube complexes.
{Proc. London Math. Soc. (3)}, 71(3) (1995) 585--617.


\bibitem{W}D.T.~Wise  Cubulating Small Cancellation Groups, Preprint, http://www.gidon.com/dani/tl.cgi?athe=pspapers/SmallCanCube.ps.

\bibitem{Yu1} G.~Yu, The Novikov conjecture for groups with finite asymptotic dimension. Ann. of Math. (2) 147 (1998), no. 2, 325--355.

\bibitem{Yu2} G.~Yu,  The coarse Baum-Connes conjecture for spaces which admit a uniform embedding into Hilbert space. Invent. Math. 139 (2000), no. 1, 201--240.


\end{thebibliography}
\end{document}